%% file: trinoids.tex
\documentclass[12pt,a4paper,notitlepage] {article}
\def\institute#1{}
\usepackage{amssymb}
\usepackage{stmaryrd}
\usepackage{amsthm}
\usepackage{artcl}
\def\qqed{}
\def\ack{\ackP}

\usepackage{mathe}
\usepackage{amsfonts}

\usepackage{graphicx}
\usepackage{enumerate}

\usepackage{float}

\usepackage{color}
\ifx\pdfoutput\undefined
  \usepackage
   {hyperref}
\else
  \definecolor{darkblue}{rgb}{.05,.05,.35}
  \usepackage[bookmarks, bookmarksopen, colorlinks, linkcolor=darkblue,
	pagecolor=darkblue, citecolor=darkblue, urlcolor=darkblue,
	pdfpagemode=UseOutlines
  ]
  {hyperref}
\fi

\def\M{\mathcal{M}}
\def\T{\mathcal{T}}

\begin{document}
\renewcommand{\thefootnote}{\fnsymbol{footnote}}

\author{Andreas Balser}
\title{Towards a classification of CMC-1 Trinoids in hyperbolic space via
conjugate surfaces\thanks{MSC2000:
53A10, 53C42, 53A35}}
\institute{Ludwig-Maximilians-Universit\"at M\"unchen, Mathematisches
Institut, Theresienstra\ss e 39, 80333 M\"unchen,
\email{balser@mathematik.uni-muenchen.de}} 

\maketitle

\begin{abstract}
We derive necessary conditions on the parameters of the ends of a CMC-1 trinoid
in hyperbolic 3-space~$\H^{3}$ with symmetry plane by passing to its
conjugate minimal surface. Together with 
\cite{benoitTrinoids}, 
this yields a classification of generic symmetric trinoids.
We also discuss the relation to other classification results of
trinoids in \cite{bobenkoTrinoids} and \cite{UY2}.

To obtain the result above, we show that the conjugate minimal surface
of a 
catenoidal CMC-1 end in $\H^{3}$ with symmetry
plane is
asymptotic to a 
suitable helicoid.
\end{abstract}


\section{Introduction}\label{introSection}
A minimal surface in $\R^{3}$ can be presented by its
\emph{Weierstrass data}, i.e.\ as a map $\Phi_{W}:\Sigma \rightarrow
\R^{3}$, where $\Sigma$ is a Riemann surface, and $\Phi_{W}$ depends
on~$(g,\omega) $, a meromorphic function and a holomorphic 1-form on
$\Sigma$. 

Given a minimal surface, one can consider its \emph{associate
(minimal) surface}, which is determined by the Weierstrass data
$(g,i\omega )$. 

Bryant found a representation of constant mean curvature 1 (CMC-1)
surfaces in $\H^{3}$ depending on the same data (see
\cite{bryant}). Therefore, we call 
CMC-1 surfaces in $\H^{3}$ \emph{Bryant surfaces}. A Bryant surface
has a \emph{minimal cousin}, the minimal surface determined by the
same data $(g,\omega)$. 

Given a Bryant surface, we define its \emph{conjugate (minimal)
surface} to be 
the associate minimal surface of its minimal cousin.

Under this construction, a principal geodesic (i.e.\ a geodesic which
is also a curvature line) on the Bryant surface corresponds to a
straight line on its conjugate surface.

We define $I:=(-\frac{1}{4},\infty)\backslash \{0 \}$, and introduce
helicoids $H_{\lambda}$, catenoids $C^{W}_{\lambda}$, and catenoid
cousins $C_{\lambda}$ parametrized by $\lambda \in I$, such that:
\newline
The helicoid $H_{\lambda}$ is the associate minimal surface of
$C^{W}_{\lambda}$, and $C^{W}_{\lambda}$ is the minimal cousin of the
Bryant surface $C_{\lambda}$.

It is known that an end of a Bryant surface is asymptotic to some
catenoid cousin or to a horosphere 
(\cite[Thm.\ 10]{geometryBryantSurfaces}).

Clearly, one would expect that the conjugate surface of a catenoidal
Bryant end is asymptotic to a suitable helicoid. However, this is not
immediate, since the Bryant cousin relation is given by a second-order
description only. For the similar situation of relating CMC-1 surfaces in
$\R^{3}$ to minimal surfaces in $S^{3}$, there exists a first-order
description. Using this, it is possible to conclude that asymptotics
is preserved in this case (see \cite{triund}).

For our situation, we show in section \ref{catenoidalEnds}
that if a catenoidal end has a symmetry plane,
then the asymptotics is indeed preserved:

\begin{theorem}\label{symmetricEnds}
Let $E'$ be a symmetric Bryant end asymptotic to $C_{\lambda}$ for some
$\lambda \in I$. Then the conjugate minimal surface $E'^{c}$ is
asymptotic to $H_{\lambda}$. 
\end{theorem}

In section \ref{trinoidSection}, we turn our attention to CMC-1
\emph{trinoids} in $\H^{3}$.  
I.e., we examine Bryant surfaces of genus zero with three ends, all of
which are catenoidal. \\
We study \emph{symmetric} trinoids, i.e.\ trinoids which have a
symmetry plane (determined by the 
asymptotic boundary 
points of their ends). 

It follows from the classification by \cite{UY2} that every (generic)
trinoid is  
symmetric; since we present a different approach to this
moduli problem, we do not use this result.
One should look for a direct geometric proof that every properly
immersed CMC-1 surface in $\H^{3}$ of genus zero and three ends has a
symmetry plane.

\smallbreak

A symmetric trinoid can be cut open along its symmetry plane
to obtain two simply connected pieces. The conjugate surface of such a
piece is a minimal surface bounded by three lines.
Surfaces of this kind were already examined by Riemann (see
\cite[sec.\ 17]{Rie} or \cite{darboux}).

Using Theorem \ref{symmetricEnds}, this yields a
necessary condition on the parameters of a 
generic trinoid:

Let $J:=(0,\infty )\backslash \{\pi \}$; for a real number $\varphi$,
we call 
\[
r(\varphi):=\min_{n\in\Z}|\varphi +2n\pi|
\]
the 
\emph{reduced angle}\index{reduced angle}\index{r@$r(\varphi )$} of
$\varphi$.  
Furthermore, let $\T$ be the set of interior points of the tetrahedron
with vertices  
$(\pi , 0,0)$, $(0,\pi ,0)$, $(0,0,\pi )$, $(\pi ,\pi ,\pi)$. Then we
have: 

\begin{theorem}\label{necCond}
If there exists a symmetric trinoid corresponding to the parameter triple 
$(\varphi_{1},\varphi_{2},\varphi_{3})\in (J\backslash \pi \Z )^{3}$,
then for the triple of reduced angles holds (in the generic case):
\[
(r(\varphi_{1}),r(\varphi_{2}),r(\varphi_{3}))\in \T.
\]
\end{theorem}

On the other hand, minimal surfaces bounded by three lines are
constructed in \cite{benoitTrinoids}.
His main result is:

\begin{theorem}[{\cite[Thm.\ 49]{benoitTrinoids}}]\label{danielThm}
Let $(\varphi_{1},\varphi_{2},\varphi_{3})\in (J\backslash \pi \Z
)^{3}$, and assume that
$(r(\varphi_{1}),r(\varphi_{2}),r(\varphi_{3}))$ lies in $\T$. 

Under a certain polynomial condition (in the $\varphi_{i}$), there is
a corresponding symmetric trinoid which arises from a minimal disk
bounded by three lines.
\end{theorem}

In section \ref{compareSection}, we compare the conditions given by
the theorems above to the conditions found in \cite{bobenkoTrinoids}
and \cite{UY2},
and find that they are essentially the same:

\begin{corollary}\label{equivalentBobenko}
The conditions of Theorems  \ref{necCond} and \ref{danielThm} are
equivalent to those given by \cite{bobenkoTrinoids}.
\end{corollary}


For symmetric parameter triples $(\varphi ,\varphi ,\varphi )$ with
$\varphi \in (\pi /3,\pi )$, one can construct the minimal surface
using a sequence of Plateau solutions, and show that it corresponds to
a trinoid; for details, see
\cite{balserDiplom}.




\def\ackII{I would like to thank Karsten Gro\ss e-Brauckmann,
who supervised my diploma thesis \cite{balserDiplom}; it was the
starting point for this article. 
}

\def\ackIII{Furthermore, I would like to thank Beno\^\i t Daniel for
useful discussions about his recent article \cite{benoitTrinoids} and
my diploma thesis.}

\def\ackP{

\paragraph{Acknowledgements.}

\ackII

\ackIII}

\def\ackM{
\begin{acknowledgement}

\ackII

\ackIII
\end{acknowledgement}
}

\ack

\section{Preliminaries}
\renewcommand{\thefootnote}{\arabic{footnote}}

In this section, we present the material from the beginning of section
\ref{introSection} in more detail.

First, we recall the Weierstrass representation for minimal surfaces
in $\R^{3}$ and the Bryant representation for CMC-1 surfaces in
$\H^{3}$:

The Weierstrass representation Theorem says that every minimal surface can be
conformally parametrized as 
\begin{align}\label{weierstrass}
\Phi_{W}(z)=\Re \oint_{z_{0}}^{z}\Kr{(1-g^{2})\omega
,i(1+g^{2})\omega, 2g\omega } ,
\end{align}
where $z$ is in $\Sigma$,  the parametrizing Riemann surface (possibly with
boundary), and $g$ (resp.\ $\omega$) is a
meromorphic function (resp.\ a holomorphic
1-form)  on~$\Sigma$. Furthermore, $g$ has a pole  of order $k$ in $z$  if
and only if $\omega$ has
 a zero of order~$2k$ in $z$.

The function $g$ has a geometric meaning: it is the stereographic
projection of the Gauss (or normal) map of the minimal surface
$\Phi_{W}$.

The pair $(g,\omega  )$ is called the \emph{Weierstrass data} of
$\Phi_{W}$. 

Conversely, Weierstrass data on a Riemann surface $\Sigma$ defines a
minimal immersion from the universal cover $\tilde{\Sigma}$ into
$\R^{3}$ via
\eqref{weierstrass}. 

To a minimal surface $\Phi_{W}$ with Weierstrass data $(g,\omega )$,
one has its \emph{associate surface} $\bar{\Phi}_{W}$, which is given
by the Weierstrass data $(g,i\omega )$; it turns out that $\Phi_{W}$
and $\bar{\Phi}_{W}$ are (locally) isometric. Note that if~$\Phi_{W}$ is
defined on~$\Sigma$, it may happen that $\bar{\Phi}_{W}$ is defined on
$\tilde{\Sigma}$ only.

For more details on the Weierstrass representation, we refer the reader to
\cite[\S 8]{ossermannSurvey}.
\medbreak

In his seminal paper \cite{bryant}, Bryant showed that there is a
representation of CMC-1 surfaces in
$\H^{3}$ using exactly the same data as the Weierstrass
representation. Thus, to 
a minimal surface $\Phi_{W}$ with Weierstrass data $(g,\omega )$, one
obtains a \emph{Bryant cousin} $\Phi_{B}$; vice versa, every CMC-1
surface in $\H^{3}$ has a \emph{minimal cousin}. 
The surfaces $\Phi_{W}$ and $\Phi_{B}$ are (locally) isometric, and
their Gauss maps agree.

\begin{definition}
Following Rosenberg, we define a \emph{Bryant surface}\index{Bryant
surface} to be an 
immersed CMC-1 surface in 
$\H^{3}$.
\end{definition}

\begin{definition}
Given a simply connected Bryant surface $M$, 
we define its 
\emph{conjugate surface}\index{conjugate surface} $M ^{c}$ 
as the associate surface of $M $'s minimal cousin (in $\R^{3}$).
\end{definition}

Since the Bryant
relation is a special case of Lawson's correspondence (\cite{UYlawson})%
, a principal geodesic (i.e.\ a geodesic which is also a curvature
line) corresponds to a principal 
geodesic under the Bryant cousin relation. Under the associate
construction, principal geodesics 
go to straight lines. Thus principal geodesics on $M$
are mapped to straight lines by $M^{c}$.
\bigbreak

\begin{example}\label{catCousExample}
We introduce our notation for
the helicoids, the catenoids, and
the catenoid cousins:

Para\-met\-rize the surfaces by $\Sigma =\C$:
For $0\not =\lambda\in \R$, the catenoid
$C_{\lambda}^{W}$ is the minimal surface with Weierstrass data $g=\exp 
(z)$, $\omega =\lambda \exp (-z)dz$, and the helicoid $H_{\lambda}$ is
its associate surface, with Weierstrass data $g=\exp 
(z)$, $\omega =\lambda i\exp (-z)dz$.

The formula for $H_{\lambda}$ is 
\begin{align}\label{helicoidFormula}
H_{\lambda }(x+iy)=2\lambda \Kr{\begin{array}{c}
\sinh x\sin y\\
-\sinh x\cos y\\
-y
\end{array}}.
\end{align}
If $\lambda \in I$, where $I:=(-\frac{1}{4},\infty )\backslash \{0
\}$, we call the Bryant cousin of 
$C_{\lambda}^{W}$ a 
\emph{Catenoid Cousin} $C_{\lambda}$.
Formulas for catenoid cousins $C_{\lambda}$ in the upper halfspace
model ($\H^{3}=\{(u+iv,w)\st u,v\in \R , w>0 \}$) are given 
in \cite[sec.\ 11]{rosen}: The surfaces are again parametrized by $\C$;
every line with constant imaginary part parametrizes a principal
geodesic from the end of $C_{\lambda}$  at $0$ to the end at~$\infty$
in~$\H^{3}$. Set 
$a:=\sqrt{1+4\lambda}$; 
then the formula for 
$C_{\lambda}(x+iy)$ is given by
\begin{align*}
\begin{array}{r@{~}l}
u+iv&=\dfrac{ -\lambda (e^{x}+e^{-x})e^{ax}}{\Kr{\frac{1}{2}+\lambda
-\frac{1}{2}a }e^{-x}+\Kr{\frac{1}{2}+\lambda
+\frac{1}{2}a}e^{x}}e^{iay} \\[4mm]
w&=\dfrac{ae^{ax}}{\Kr{\frac{1}{2}+\lambda
-\frac{1}{2}a}e^{-x}+\Kr{\frac{1}{2}+\lambda +\frac{1}{2}a}e^{x}}
\end{array}
\end{align*}
Note that the parametrization of $C_{\lambda}$ is
periodic with period $\frac{2\pi i}{\sqrt{1+4\lambda}}$.
\end{example}

Let $J:=(0,\infty )\backslash \{\pi  \}$\index{J@$J$}, and define the
bijective function 
$\tilde{\varphi} :I\rightarrow J$\index{4@$\tilde{\varphi}(\lambda )$}
 by $\tilde{\varphi} (\lambda ):=\frac{\pi 
}{\sqrt{1+4\lambda }}$.

We remark that the Catenoids (and Catenoid Cousins)
$C_{\lambda}^{(W)}$ can alternatively
be described by the Weierstrass data $g=z^{\alpha }, \omega
=\frac{1-\alpha ^{2}}{4\alpha }z^{-1-\alpha }$ on~$\C^{*}$, where
$\pi \alpha  =\tilde{\varphi}(\lambda )$; see \cite[Ex.\ 1.5]{Earp}.

\section{Symmetric catenoidal Bryant ends and their conjugate surfaces}
\label{catenoidalEnds}
In this section, we show that the conjugate minimal surface of a
catenoidal Bryant end with a symmetry plane is asymptotic to the
corresponding helicoid.

\begin{definition}
An \emph{annular Bryant end} is a Bryant surface with domain
$\{0<|z|\leq 1\}$ (or equivalently any other punctured disk with boundary).
\end{definition}

Recall that  a properly embedded Bryant
annular end in $\H^{3}$ is asymptotic to some catenoid cousin or to a
horosphere \cite[Thm.\ 10]{geometryBryantSurfaces}.

\begin{definition}
An annular Bryant   end  is called
\emph{catenoidal}\index{catenoidal end} if it 
is properly embedded and asymptotic to a catenoid cousin.

An  annular Bryant end is called \emph{symmetric} if
it is  properly embedded and has a symmetry plane.
\end{definition}
\begin{definition}
A \emph{minimal end bounded by rays} is a properly immersed mini\-mal
surface in $\R^{3}$ 
with domain $\{0<|z|\leq 1, \Im z\geq 0\}$, such that $[-1,0)$ and
$(0,1]$ are mapped to (monotonically parametrized) rays.
\end{definition}

For a real number $\varphi$, we call
$r(\varphi):=\min_{n\in\Z}|\varphi +2n\pi|$ the 
\emph{reduced angle}\index{reduced angle}\index{r@$r(\varphi )$} of
$\varphi$.  

The following Lemma is a slight generalization of 
\cite[L.\ 7]{benoitTrinoids}: 
\begin{lemma}\label{benoitLemma}
Let $X$ be a minimal end bounded by horizontal rays with vertical
limit normal for $z\rightarrow 0$. Assume that $X$ is contained in a vertical
slab (i.e.\ the vertical component of $X$ is bounded), that the
stereographic projection $g$ of the Gauss map of $X$ 
satisfies $g\sim 
z^{\alpha}$ for $z\rightarrow 0$ with $0<\alpha \not =1$, and that the
vertical components of the boundary rays are
$|\tilde{\varphi}^{-1}(\pi \alpha )|\pi \alpha$ apart. 
 Then~$X$ is
asymptotic to (part of) $IH_{\lambda}$ for $\lambda =\tilde{\varphi}^{-1}(\pi
\alpha )$ and some orientation preserving isometry~$I$ of $\R^{3}$. 
\end{lemma}
\begin{proof}
First we note that we can conclude from the proof of 
\cite[L.~7]{benoitTrinoids} that the angle between the boundary
rays is the reduced angle $r(\pi \alpha )$.
Additionally, observe that the (vertical) distance of the boundary
rays is by assumption the
distance of two lines in $H_{\lambda}$, where one has to be rotated by
angle $\pi \alpha$ \emph{in} $H_{\lambda}$ to be mapped to the other
one (cf.\ formula \eqref{helicoidFormula}).

If the boundary rays are not parallel, the claim is just
\cite[L.\ 7]{benoitTrinoids}. 
The case of parallel boundary rays is not covered there;  
 however, its proof still works in this
case by our assumptions on the limit normal, $X$ being contained in a
slab, and the  vertical 
distance of the rays.
\qqed 
\end{proof}

\begin{proof}[Proof of Theorem \ref{symmetricEnds}]
It suffices to consider one symmetric piece $E$ of $E'$ bounded by principal
geodesics (the curves of intersection with the symmetry plane). Let
$E^{c}$ denote the conjugate surface of this half. We 
assume~$E^{c}$ to be parametrized by $D:=\{0<|z|\leq 1, 
\Im z\geq 0\}$. 
By \cite{geometryBryantSurfaces}, $E'$ has a well-defined limit
normal, which we may assume to be vertical.

Then $E^{c}$ also has a vertical limit normal, so it is a minimal end
bounded by horizontal rays.
We show that $E^{c}$ is contained in a vertical slab:

By \cite{Earp}, we may assume the Weierstrass data of $E'$ to be of  the form 
\[
g = z^{\alpha}(g_{0}+zg_{1}(z)),\quad \omega
=z^{-1-\alpha}(w_{0}+zw_{1}(z)) 
\]
with $g_{0},w_{0}\in \C$ such that
$g_{0}w_{0}=\frac{1-\alpha^{2}}{4\alpha}$, and 
holomorphic functions $g_{1},w_{1}$ on $\{|z|\leq 
1\}$ (where $\pi \alpha  =\tilde{\varphi}(\lambda)$, in particular
$0<\alpha  \not =1$). 

Choose $z_{0}\in (0,1]\subset D$; the third component of $E^{c}$ is
the negative of the imaginary part of the following integral:
\begin{align*}
\oint_{z_{0}}^{z}2g\omega &=2\oint_{z_{0}}^{z}
\xi^{-1}\Bigl(g_{0}w_{0} + \xi
\bigl( g_{0}w_{1}(\xi ) +w_{0}g_{1}(\xi )+\xi w_{1}(\xi)g_{1}(\xi)
\bigr)\Bigr)d\xi \\ 
&=\frac{1-\alpha^{2}}{2\alpha}\oint_{z_{0}}^{z}\xi^{-1}d\xi
+2\underbrace{\oint_{z_{0}}^{z} g_{0}w_{1}(\xi ) +w_{0}g_{1}(\xi
)+\xi w_{1}(\xi)g_{1}(\xi)d\xi}_{=:\,C(z)}  
\end{align*}
Hence, $E^{c}$ is contained in a vertical slab, since $C$ is bounded
on $D$ and the first summand corresponds to the third component of
$H_{\lambda}$. 
Observe that the imaginary part of the first summand above is 0 
for $z\in (0,1]$ and constant for $z\in [-1,0)$.
We show that $\Im C(z)=0$ for $z\in [-1,0)\cup (0,1]$: This is
clear for $z\in (0,1]$, since $z_{0}\in (0,1]$, and a horizontal
ray is  parametrized. Similarly, $\Im C(z)\equiv C_{2}$ for $z\in
[-1,0)$ since this 
parametrizes another horizontal ray. 
Thus $\Im\oint_{1/n}^{-1/n}g_{0}w_{1}(\xi )
+w_{0}g_{1}(\xi )+\xi w_{1}(\xi)g_{1}(\xi)d\xi$ is constant
(i.e.\ independent from $n$ and the path in $D$ from $\frac{1}{n}$ to
$-\frac{1}{n}$) and we have 
\[
C_{2}=\Im \lim_{n\rightarrow \infty}\oint_{1/n}^{-1/n}g_{0}w_{1}(\xi )
+w_{0}g_{1}(\xi )+\xi w_{1}(\xi)g_{1}(\xi)d\xi = 0.
\]

This shows that the two boundary rays of $E^{c}$ have positive
vertical distance, which is equal to the distance of
corresponding lines on $H_{\lambda}$. Now the conclusion follows via
Lemma \ref{benoitLemma}.
\qqed
\end{proof}

\begin{corollary}\label{boundaryLines}
Let $E'$ be a symmetric Bryant end which is asymptotic to
$C_{\lambda}$, and
let $E$ be a symmetric  piece of $E'$ as above.
If $\varphi :=\tilde{\varphi}(\lambda)\not \in \pi \Z$, we have: The
boundary rays 
$l_{1},l_{2}$ of $E^{c}$
are contained in 
$IH_{\lambda}$ for some orientation-preserving isometry
$I$ of $\R^{3}$.
In particular, the angle between the ends of $l_{1}$ and
$l_{2}$ 
 is 
$r(\varphi )$. The
distance of these two lines is 
$h(\varphi ):=|\lambda|\varphi$.
\end{corollary}

\begin{proof}
First we note that $h:J\rightarrow \R $ 
is
well-defined, because $\tilde{\varphi} $ is a bijective function
(in fact $h(\varphi )=|\frac{\pi^{2}}{4\varphi}-\frac{\varphi}{4}|$).

The claim follows immediately from the proof of Theorem \ref{symmetricEnds}
and the formulas for helicoids.
\qqed\end{proof}

In case $\tilde{\varphi}(\lambda)\in \pi \Z$, the boundary rays are
parallel, and we have a
 lower bound on their distance (by the distance of parallel lines in
the corresponding helicoid).


\section{Trinoids}\label{trinoidSection}

\begin{definition}\label{trinoidDef}
We define a \emph{trinoid}\index{trinoid} to be a properly immersed
Bryant  surface 
of genus zero with three ends, all of which are 
ca\-te\-no\-idal. 
A \emph{symmetric trinoid} is a trinoid $T$ 
which
has a symmetry plane $P$ such that the asymptotic endpoints of $T$ are
contained in the asymptotic boundary of $P$. 


Denote by $\M$ the space of symmetric trinoids with ends marked 
by 1, 2, 3, up to isometry (respecting the marks of the ends).
\end{definition}

Observe that the symmetry plane $P$ is uniquely determined if the
asymptotic endpoints are distinct.

Pictures of trinoids can be found at \newline
\texttt{http://www-sfb288.math.tu-berlin.de/\~{}bobenko/Trinoid/webimages\newline
.html};  
see also \cite{bobenkoTrinoids}.

\begin{definition}
We can define the map $\Psi : \M\rightarrow J^{3}$ sending a
trinoid to the triple $(\varphi _{1},\varphi  _{2},\varphi _{3})\in
J^{3}$, where $\varphi _{i}=\tilde{\varphi} (\lambda _{i})$, and $\lambda _{i}$
is the parameter of end~$i$.
\end{definition}

\begin{lemma}
Any properly embedded Bryant surface $M$ of genus zero with three ends is
a symmetric trinoid.
\end{lemma}
\begin{proof}
By Theorem \cite[Thm.\ 12]{geometryBryantSurfaces}, every end is
catenoidal, and by  
\cite[Thm.\ 11]{geometryBryantSurfaces}, the three asymptotic boundary
points are distinct
and~$M$ is a bigraph over the plane containing them. \qqed
\end{proof}

We expect that the Lemma above generalizes to Alexandrov-embedded
Bryant surfaces.



Note that a trinoid is a map $S^{2}\backslash \{x_{1},x_{2},x_{3}
\}\rightarrow \H^{3}$, where $x_{1},x_{2},x_{3}$ are distinct and
correspond to the ends 1, 2, 3 respectively.

\begin{lemma}\label{l12l23l31}
For a symmetric trinoid $M\in \M$, there is a unique principal geodesic of~$M$
joining $x_{1}$ to $x_{2}$, which we denote by $l_{12}$. Similarly,
there is a unique principal 
geodesic $l_{23}$ joining $(x_{2},x_{3})$ and a unique principal geodesic
$l_{31}$ joining 
$(x_{1},x_{2})$. \\
Considering  $l_{12},l_{23},l_{31}$ as subsets of
$S^{2}$, we have that $S^{2}\backslash (l_{12}\cup l_{23}\cup l_{31}\cup
\{x_{1},x_{2},x_{3} \})$  consists of exactly two components.
\end{lemma}
\begin{proof}
It is  known that a principal geodesic is contained in a plane of symmetry
 of $M$ (cf.\ \cite[Prop.\ 3.2]{Earp}). 
We conclude that the three lines we are looking for need to be
contained in $P$, the symmetry plane of $M$ from the definition.

Consider the graph $G$ in $S^{2}$ with vertices $V:=\{x_{1},x_{2},x_{3} \}$,
and edges the principal geodesics of $M$ contained in $P$ which start
or end in $V$ (observe that both asymptotic ends of such a principal
geodesic are in $V$).

Since every end of $M$ is embedded, every vertex has degree two. 
Edges cannot intersect: Tangential contact is excluded by uniqueness
of geodesics, and transversal intersection is impossible since $M$
intersects $P$ orthogonally near every point of $G$.

Thus, $G$ consists of one, two, or three loops in $S^{2}$. Reflection
in $P$ maps every component of $S^{2}\backslash G$ to an other
component. Since all elements of $V$ are fixed points of this
reflection, $G$ consists of one loop only.
\qqed\end{proof}

\begin{corollary}\label{halfspaceNearPrinc}
Consider a symmetric trinoid $M$ and its symmetry plane~$P$.
Then there is a neighborhood $N$ of $l_{12}\cup l_{23}\cup l_{31}\cup
\{x_{1},x_{2},x_{3} \}$ in $S^{2}$ such that $M(N)\cap P=M(l_{12}\cup
l_{23}\cup l_{31})$. In particular: Near its boundary, each component of 
$S^{2}\backslash (l_{12}\cup l_{23}\cup l_{31}\cup
\{x_{1},x_{2},x_{3} \})$ is mapped to a component of~$\H^{3}\backslash
P$.

If $M$ is embedded, each component of $S^{2}\backslash (l_{12}\cup
l_{23}\cup l_{31}\cup 
\{x_{1},x_{2},x_{3} \})$ is mapped into a halfspace of $\H^{3}\backslash P$.
\qed\end{corollary}


Given a symmetric trinoid $M$, we can (by an orientation-preserving isometry)
assume that its symmetry 
plane is the equatorial plane $E=\{x_{3}=0 \}$ of the Poincar\'e disk
model (lying inside $\R^{3}$). Further, we can 
assume that the ends are marked increasingly if one looks from above
(i.e.\  the direction of positive $x_{3}$).

\begin{definition}
Given  a symmetric trinoid $M$, we divide its domain $S^{2}\backslash
\{x_{1},$ $x_{2},x_{3} \}$ into two components 
along $l_{12},l_{23},l_{31}$
, and we
define $M^{+}$\index{M@$M^{+}$} to be the restriction of $M$ to the
closure of the component which is 
mapped to the upper half space near $l_{12},l_{23},l_{31}$,
 if $M$ is put in the Poincar\'e
model in the 
way explained above.
\end{definition}

\begin{figure}
\centerline{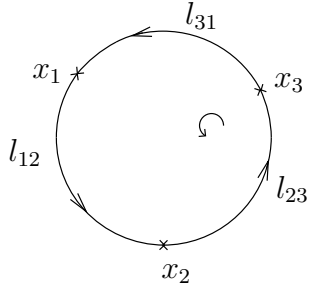}
\caption{The domain of $M^{+}$ and $M^{c}$.}\label{domain}
\end{figure}

So $M^{+}$ is a map 
$\bar{D}\backslash \{x_{1},x_{2},x_{3} \}\rightarrow \H ^{3}$, where
 $x_{1}$, $x_{2}$,
$x_{3}$ are distinct points in $\partial D$ (and $D$ is the closed
unit disk).\\
We choose the orientation on $D$ and its
boundary as depicted in Figure \ref{domain}.

Note that $M^{+}$ is well-defined up to orientation-preserving
hyperbolic isometries leaving the upper half-space in the Poincar\'e
disk model  invariant.

Define $M^{c}:=(M^{+})^{c}$ to be the conjugate minimal surface of $M^{+}$.

Since principal geodesics on a Bryant surface correspond to straight
lines on its conjugate minimal surface, we have:
\begin{lemma}
Let $M$ be  a symmetric trinoid. Then $M^{c}$ is a minimal surface
bounded by three straight lines.\qed
\end{lemma}

We mention another interesting fact about $M^{+}$:
\begin{proposition}
If $M^{+}$ is embedded, then so is $M$.
\end{proposition}
\begin{proof}
Assume the symmetry plane of $M$ to be the equatorial plane $E$ as
before. 
We apply the Alexandrov-reflection technique (see, for example,
\cite{constantMeanHypersurfaces}):
Using a (continuous) family of planes which foliate the upper
half-space, we conclude that the normal of $M^{+}$ at any point 
$p\in M^{+}\cap E$ has  non-positive vertical coordinate.

Similarly, we use a family of planes foliating the lower half-space to
find that for a point $p\in M^{+}\cap E$, the normal has
to have non-negative vertical coordinate.

Thus, every component of $M\cap E$ is a principal geodesic, i.e.\ a
curve of planar reflection (by \cite[Prop.\ 3.2]{Earp}).
So $M^{+}$ is cut off wherever it reaches~$E$ (observe that there are
no closed principal geodesics since $M$ has genus zero), and it does not
 intersect the lower half-space; so $M$ is embedded.\linebreak
\vspace*{1in}
 \qqed\end{proof}


\section[Necessary conditions]{Necessary conditions on the
constellation of boundary lines}

In this section, we use the information about the constellation of
lines which bound $M^{c}$ to obtain a necessary condition on the
parameter triple $\Psi (M)$ in the generic case.


\begin{definition}\label{admissDef}
A triple $(\lambda_{1},\lambda_{2},\lambda_{3})\in I^{3}$ is called a
\emph{parameter triple}\index{parameter triple}. 

A triple of oriented lines $(l_{12},l_{23},l_{31})$ in $\R^{3}$ is called
\emph{admissible constellation}\index{admissible constellation} if 
\begin{enumerate}[(i)]
\item 
There exists a parameter triple
$(\lambda_{1},\lambda_{2},\lambda_{3})\in I^{3}$ and
orientation-pre\-ser\-ving isometries $I_{1},I_{2},I_{3}$ of $\R^{3}$ such that
\[
\{l_{12},l_{31} \}\subset I_{1}(H_{\lambda_{1}}),\quad  
\{l_{23}, l_{12}\}\subset I_{2}(H_{\lambda_{2}}),\quad 
\text{and } \{ l_{31}, l_{23} \}\subset
I_{3}H_{\lambda_{3}}.\]
\item 
Rotating $l_{i(i+1)}$
inside $I_{i}(H_{\lambda_{i}})$ maps $l_{i(i+1)}$ to $l_{(i+2)i}$ with the
opposite orientation (for $i\in \Z_{3}$).
\item The distance of $l_{i(i+1)}$ and $l_{(i+2)i}$ is 
$h\circ \tilde{\varphi}(\lambda_{i})$ (for $i\in \Z_{3}$).
\end{enumerate}

A triple $(\varphi _{1},\varphi _{2},\varphi _{3})\in J^{3}$ of
angles is 
called \emph{admissible}, if there exists an admissible constellation
with parameter triple $(\tilde{\varphi}^{-1}(\varphi_{1}),
\tilde{\varphi}^{-1}(\varphi_{2}),\tilde{\varphi}^{-1}(\varphi_{3}))$.

An admissible triple is called \emph{generic}
 if there is a 
corresponding admissible constellation such that the lines
are not contained in parallel
planes.
The triple is called \emph{parallel}\index{parallel triple} otherwise.

\begin{remark}
Note that a general triple of three oriented lines is determined (up
to the action of $SO(3)$) by the oriented distances and the angles. We
have the  restriction that the  distance and
angle match, i.e.\ every pair of lines can be put into a suitable
helicoid.  
\end{remark}

We define $\T$ to be the set of interior points of the tetrahedron
with vertices  
$(\pi , 0,0)$, $(0,\pi ,0)$, $(0,0,\pi )$, $(\pi ,\pi ,\pi)$.
\end{definition}


Sketches of admissible constellations can be found in \cite{balserDiplom}.

From Corollary \ref{boundaryLines}, we have:
\begin{lemma}\label{trinoidsAdmissConst}
For any symmetric trinoid $M\in
\M$ with $\tilde{\varphi}(\lambda_{i})\not \in \pi \Z$ for $i\in \{1,2,3 \}$,
the triple $\Psi  (M)$ is admissible.\qed 
\end{lemma}

\begin{theorem}\label{admiss}
A triple $(\varphi_{1},\varphi_{2},\varphi_{3})\in J^{3}$ is
a generic admissible triple if and only if
$(r(\varphi_{1}),r(\varphi_{2}),r(\varphi_{3}))\in \T$. For every
triple of that kind,  there are exactly two generic admissible
constellations of lines 
in $\R^{3}$ (modulo $SO(3)$).
\end{theorem}

Proofs can be found in
\cite[Prop.\ 9]{benoitTrinoids}; and \cite{balserDiplom}.
Essentially, the conditions on the angles correspond to the condition
that the directions of the lines form a spherical triangle (after
identifying the unit tangent spheres of $\R^{3}$ via parallel translation).

\smallbreak
\begin{proof}[Proof of Theorem \ref{necCond}]
The theorem  follows from Lemma
\ref{trinoidsAdmissConst} and Theorem \ref{admiss}. 
\end{proof}


\begin{remark}
One can show that an admissible triple corresponds \emph{either} to
generic \emph{or} to parallel constellations, and that the triple of
reduced angle lies in the boundary of $\T$ in the parallel case, see
\cite{balserDiplom}. Hence the name \emph{generic} is justified.
\end{remark}


\section{Comparing to related results}\label{compareSection} 
In this section, we compare the conditions obtained by
\cite{benoitTrinoids} and our results with the results in 
in \cite{bobenkoTrinoids} and \cite{UY2}.

Consider the presentation of
catenoid cousins in \cite[Ex.\ 2]{bryant}. Bryant parametrizes
catenoid cousins with a parameter $-\frac{1}{2}<\mu_{B}\not =0$. 

\begin{lemma}
The catenoid cousin given by Bryant's parameter $\mu_{B}$ is
$C_{\lambda}$, where 
$\lambda =\tilde{\varphi}^{-1}(\pi (2\mu_{B} +1))$.
\end{lemma}
\begin{proof}
Bryant computes the total curvature  of a catenoid
cousin to be \newline
$-4\pi (2\mu_{B} +1)$. 
A standard catenoid has total curvature $-4\pi$. Since  a Bryant
surface is locally isometric to its minimal cousin, a catenoid
cousin $C_{\lambda}$ has total curvature $-4\pi \cdot
\frac{1}{\sqrt{1+4\lambda}}$ (see Example \ref{catCousExample}). The
claim follows.  
\end{proof}

Next, we  trace back the
relationship between our parameters and the parameters in
\cite{bobenkoTrinoids}.

In \cite[sec.\ 4]{bobenkoTrinoids}, catenoid cousins are parametrized
by a parameter $0<\lambda_{BPS}\not =\frac{1}{2}$. Comparing the formulas for
catenoid cousins given by Bryant and \cite{bobenkoTrinoids}, 
 we obtain
$\lambda_{BPS}=\mu_{B}+\frac{1}{2}$; hence, the catenoid cousin
described by the parameter $\lambda_{BPS}$ is $C_{\lambda}$, where 
\begin{align}\label{relation}
\lambda = \lambda (\lambda_{BPS}) =\tilde{\varphi}^{-1}(2\pi
\lambda_{BPS}).
\end{align}


They consider $|\{\lambda _{BPS,i} \}|$, where $\{\cdot \}$ stands for
the 
fractional part of a number in $[-\frac{1}{2},\frac{1}{2})$. 

The main result in \cite{bobenkoTrinoids} is: 
\begin{theorem}[{\cite[Prop.\ 2]{bobenkoTrinoids}}] \label{bobenkoThm}
For given parameters $p_{i},q_{i}$, where $i\in \{0,1,\infty \}$, in the
generic case, it is necessary for the existence of 
a trinoid that the numbers 
$\Delta_{i}:=|\{\lambda _{BPS,i}\}|$ satisfy the conditions 
\begin{align*}
\Delta_{0}+\Delta_{1}+\Delta_{\infty}&>\frac{1}{2}\\
\Delta_{0}+\Delta_{1}-\Delta_{\infty}&<\frac{1}{2}\\
\Delta_{0}-\Delta_{1}+\Delta_{\infty}&<\frac{1}{2}\\
-\Delta_{0}+\Delta_{1}+\Delta_{\infty}&<\frac{1}{2};
\end{align*}
This condition is sufficient if furthermore, certain holomorphic
spinors $P$ and $ Q$ have no common zeroes.\qed
\end{theorem}

They also show that their classification is equivalent to
\cite[Thm.~2.6]{UY2}.  

Observe that the ``generic case'' in \cite{bobenkoTrinoids} means that
the  case of half-integer $\lambda_{BPS,i}$  
is excluded (cf.\ formula (6.3), and the remark at the bottom of
page~18), so their class of generic trinoids is slightly larger than ours.

\medbreak

In \cite[Thm.\ 49]{benoitTrinoids}, the trinoids from the
classification of Umehara-Yamada (or equiv.\ Bobenko et al.) are
constructed via minimal surfaces bounded by a generic constellation of
three lines. 

\smallbreak

In view of
\eqref{relation}, we find that $|\{\lambda_{BPS} \}|$  corresponds to our 
notion of reduced angle, i.e.\ we have 
$r(\tilde{\varphi}\circ \lambda (\lambda_{BPS}))=2\pi |\{\lambda_{BPS} \}|$. 
So the necessary conditions of Theorem \ref{bobenkoThm} are the same
as those in Theorem \ref{necCond}.

Comparing our Theorem \ref{necCond}  to the main theorem of
\cite{bobenkoTrinoids}, the condition about
the common zeroes of $P,Q$ is preventing our condition from being
sufficient. 
In \cite{benoitTrinoids}, this additional condition is that his
polynomial ``$\varphi$'' of degree two 
 has no double root (for the equivalence, see \cite[proof of L.\ 16,
and page 31]{benoitTrinoids}). This condition avoids singular points on
the minimal surface and the trinoid.

\smallbreak

Hence, Corollary \ref{equivalentBobenko} follows.


\bibliographystyle{alpha}
\bibliography{dipl}

 \end{document}

%% file: 2fig1.pstex_t
\begin{picture}(0,0)%
\includegraphics{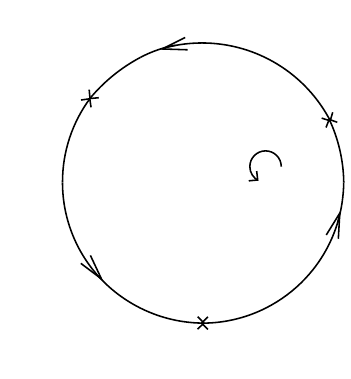}%
\end{picture}%
\setlength{\unitlength}{3947sp}%
\begingroup\makeatletter\ifx\SetFigFont\undefined%
\gdef\SetFigFont#1#2#3#4#5{%
  \reset@font\fontsize{#1}{#2pt}%
  \fontfamily{#3}\fontseries{#4}\fontshape{#5}%
  \selectfont}%
\fi\endgroup%
\begin{picture}(1679,1811)(1426,-2502)
\put(2535,-847){\makebox(0,0)[lb]{\smash{\SetFigFont{12}{14.4}{\familydefault}{\mddefault}{\updefault}{\color[rgb]{0,0,0}$l_{31}$}%
}}}
\put(3105,-1938){\makebox(0,0)[lb]{\smash{\SetFigFont{12}{14.4}{\familydefault}{\mddefault}{\updefault}{\color[rgb]{0,0,0}$l_{23}$}%
}}}
\put(2385,-2444){\makebox(0,0)[lb]{\smash{\SetFigFont{12}{14.4}{\familydefault}{\mddefault}{\updefault}{\color[rgb]{0,0,0}$x_2$}%
}}}
\put(3093,-1244){\makebox(0,0)[lb]{\smash{\SetFigFont{12}{14.4}{\familydefault}{\mddefault}{\updefault}{\color[rgb]{0,0,0}$x_3$}%
}}}
\put(1426,-1711){\makebox(0,0)[lb]{\smash{\SetFigFont{12}{14.4}{\familydefault}{\mddefault}{\updefault}{\color[rgb]{0,0,0}$l_{12}$}%
}}}
\put(1576,-1186){\makebox(0,0)[lb]{\smash{\SetFigFont{12}{14.4}{\familydefault}{\mddefault}{\updefault}{\color[rgb]{0,0,0}$x_1$}%
}}}
\end{picture}